\documentclass[11pt]{article}
\usepackage{amssymb}
\usepackage{latexsym}

\newtheorem{thm}{Theorem}[section]
\newtheorem{prop}[thm]{Proposition}
\newtheorem{lemma}[thm]{Lemme}

\newtheorem{dfn}[thm]{D\'efinition}

\begin{document}

\title{\bf  A step to Gronwall's conjecture.  } 

\author{ Jean Paul Dufour. \\}

\maketitle

\begin{abstract}
In this paper we will explore a way to prove the hundred years old Gronwall's conjecture: if two plane linear 3-webs  with non-zero curvature are locally isomorphic, then the isomorphism is a homography. 
 
Using recent results of S. I. Agafonov, we exhibit an  invariant, the {\sl  characteristic}, attached to each generic  point of  such a web, with the following property: if a diffeomorphism interchanges two such linear webs, sending a point of the first to a point of the second which have the same  characteristic, then this diffeomomorphism is locally a homography.
\end{abstract}


\section{Introduction}

In this text we work in the real projective plane. The results are the same in the complex case. We also work in the analytic case.

A plane 3-web is a triple of 1-dimensional foliations, two by two tranversal, on an open domain of the plane. In the sequel we will work only with these 3-webs : so we forgot the word {\sl plane}. Such a  3-web is called {\sl linear} if the  leaves of the foliations are rectilinear.

Let $W$ and $\bar W$ be two 3-webs, linear or not, defined respectively on domains $U$ and $\bar U.$ We say that they are isomorphic if there is a diffeomorphism from $U$ to $\bar U$ which maps every leaf of the foliations of $W$ on a leaf of the foliations of $\bar W.$

Near any point of the domain of a 3-web $W,$ there are coordinates $(x,y)$ (in general non affine) such that the foliations  are given by the verticals $ x={ constant},$ the horizontals $y={ constant}$ and  the level sets of some function $f$ (the sets of $(x,y)$ such that $ f(x,y)
={constant}).$ Such a 3-web is denoted $(x,y,f).$ In other words we can say that $W$ is locally isomorphic to some $(x,y,f).$

Attached to each 3-web there is a 2-form called its (Blaschke) {\sl curvature}. For $(x,y,f)$ this curvature is
$$\partial_x \partial_y(\log{\partial_x f\over \partial_y f}) dx\wedge dy.$$

We say that a 3-web is {\sl non flat} if its curvature doesn't vanish at any point of its domain.

\noindent{\bf  Gronwall's conjecture : } {\sl If two non flat linear 3-webs are isomorphic, then the diffeomorphism which realizes this isomorphism is a homography near any point.}
 
Probably  the recent paper of S.I. Agafonov (\cite{SA})  contains the best historical references on this subject. We will also use some notions appearing in this paper and,  particularly,  the following.

    We choose affine coordinates $(x,y)$ on an open subset $U$ of the plane. Let $W$ a linear 3-web defined on $U$ by the slopes $P,$ $Q$ and $R$ of the different foliations ($P(x,y)$ is the slope of the first foliation at the point of coordinates $(x,y)$ ...). We define following quantities
$$\Pi=(P-Q)(Q-R)(R-P),$$
$$\Delta=(P-Q)R_y+(Q-R)P_y+(R-P)Q_y,$$
with the convention:
if $f$ is a function of $(x,y)$ we denote $f_y$ its derivative with respect to $y,$ $f_{yy}$ its second derivative with respect to $y$.... We assume now that $W$ is non flat. The {\bf Lemma 1} of \cite{SA} implies that we can assume also that $\Delta$ is everywhere non vanishing. 

\begin{dfn} The characteristic of $W$ is $$car_W:=\Pi . (P_{yy}+Q_{yy}+R_{yy})/\Delta^2.$$  \end{dfn}

In \cite{SA} we can find a complete set of projective invariants (invariant up to homographies) and 
$car_W$ is the sum of three of them.

Because $ car_W$ is a projective invariant,
if a homography $\psi$ interchanges  linear 3-webs $W$ and $\bar W$,  we have 
$$car_W(M)=car_{\bar W}(\psi(M)),$$ for any point $M.$

   Our central result is the following. 
\begin{thm}\label{resultatA} 
We consider two non flat linear 3-webs $W$ and $\bar W.$ We assume that there is an isomorphism $\phi$ from the first one  to the second. If there is a point $M$ such that 
\begin{equation}\label{(1)}car_W(M)=car_{\bar W}(\phi(M))\end{equation}
 then the isomorphism is a homography near $M.$
\end{thm}

 This theorem says that Gronwall's conjecture is true  if we can always find a point $M$ satisfying condition (\ref{(1)}). We can hope to find such $M$ by a fix point method.

\section{Description of linear 3-webs near a point.}\label{descro}

 We consider a linear 3-web $W,$ and $M$ a point of its domain. Then we can prolongate the three foliations to obtain three  families of lines $\cal A,$ $\cal B$ and $\cal C.$ In the projective plane $\cal A,$ $\cal B$ and $\cal C$  envelop three curves, which may degenerate into a point. We denote respectively $A_M ,$ $B_M$ and $C_M$ the focal points on the three line passing by $M$ (the points where these lines touch the envelops). We impose now that $W$ is non flat. Then  the {\bf Lemma 1} of \cite{SA} implies also
that $(M,A_M,B_M,C_M)$ is a projective frame. Up to a homography, we can choose coordinates $(u,v)$ such that
$$M=(0,0),\ A_M=(1,1), \ B_M=(1,-1),\ C_M=(2,0).$$
Near the origin each leaf of our foliations is transversal to the $v$-axis. So our 3-web can be described as follows. 

There are three 1-variable functions  $a:x\mapsto a(x)$ (resp. $b:y\mapsto b(y)$, $c:z\mapsto c(z)$) such $\cal A$ (resp $\cal B$,  $\cal C$) consists of the lines $v=a(x)u+x$ (resp.  $v=b(y)u+y$, $v=c(z)u+z$) where $x$ (resp. $y$, $z$) is a parameter varying near  the origin.

Be aware that the functions $a$, $b$ and $c$ we just defined aren't the $a$, $b$ and $c$ of \cite{SA}.

Because of the choice of $(u,v)$, Taylor expansions of $a,$ $b$ and $c$ have the shapes
$$a(t)=1-t+a_2t^2+\cdots +a_it^i+\cdots ,$$  $$b(t)= -1-t+b_2t^2+\cdots +b_it^i+\cdots,$$ $$ c(t)=- t/2 +c_2t^2+\cdots+c_it^i+\cdots .$$

\begin{lemma}\label{car(0,0)} The following formula holds:
$$car_W(0,0)=4(a_2+b_2+c_2).$$
Moreover  $a_2+b_2+c_2$ vanishes  if and only if the curvature of $W$ at the origin  vanishes.\end{lemma}

To prove this lemma we begin to compute the 2-order   Taylor expansion of the function $x(u,v)$ (resp. $y(u,v),$  resp. $z(u,v)$) given by the implicit relation $v=a(x(u,v))u+x(u,v)$ (resp. $v=b(y(u,v))u+y(u,v),$ resp. $v=c(z(u,v))u+z(u,v)$). We find they are $-u+v-u^2+uv$ for $x(u,v),$  $u+v+u^2+uv$ for $y(u,v)$ and $v+uv/2$ for $z(u,v)$).

 Then the slope functions are $P(u,v)=a(x(u,v)),$  $Q(u,v)=b(y(u,v))$ and $R(u,v)=c(z(u,v)).$ So their 2-order Taylor expansions are respectively $$1+u-v+(1+a_2)u^2+(-1-2a_2)vu+a_2v^2,$$ $$-1-u-v+(-1+b_2)u^2+(-1+2b_2)vu+b_2v^2$$ and $$-v/2-uv/4+c_2v^2.$$

So the values of $\Pi,$ $\Delta$, $P_{yy}$, $Q_{yy}$ and $R_{yy}$ at the origin are respectively
 $2,$ $1,$ $2a_2,$ $2b_2$ and $2c_2.$ 
 This proves  the first assertion of the lemma.

The second can be proved, for example, using the formula of the curvature given in the introduction.

\section{Isomorphic linear webs.}

In order to prove  theorem \ref{resultatA}, we consider two non flat linear 3-webs $W$ and $\bar W,$ the first near a point $M$, the second near a point $\bar M.$ For  both we adopt a description as in the preceeding section : the first is described by the three local functions $a,$ $b$ and $c,$ the second is described by the local functions $\bar a,$ $\bar b$ and $\bar c.$ Via the lemma \ref{car(0,0)}, the relation $car_W(M)=car_{\bar W}(\bar M)$ writes as
$$a_2+b_2+c_2=\bar a_2+\bar b_2+\bar c_2,$$
with evident notations.

In the following of this section we propose a way to express the existence of an isomorphism from $W$ to $\bar W$ mapping $M$ to $\bar M.$

Let $x,$ $y,$ $z$ be three numbers such that  the three leaves $v=a(x)u+x,$ $v=b(y)u+y,$ $v=c(z)u+z$ of $W$ have a commun point. They are characterized by the relation
\begin{equation}\label{pointcommun}  det\left|\matrix{1 & 1 & 1\cr x & y & z \cr a(x) & b(y) &  c(z)  \cr}\right|=0\ .\end{equation}

This equation defines implicitly a function $z=f_W(x,y)$ on a neigborhood of the origin. This proves that $W$ is locally isomorphic to  $(x,y,f_W).$  Because the change of coordinates  $(u,v)\mapsto (x,y)$ is not a homography (in general),   $(x,y,f_W)$ has no reason to be linear.

Replacing respectively $a$ with $\bar a,$  $b$ with $\bar b,$  $c$ with $\bar c ,$ we can construct $f_{\bar W}$ such that $\bar W$ is isomorphic to the 3-web $(x,y,f_{\bar W}).$

Now $W$ is locally isomorphic to $\bar W$  by a diffeomorphism which maps $M$ to $\bar M,$ if and only if there is a local isomorphism, preserving the origin, from  $(x,y,f_W)$ to  $(x,y,f_{\bar W}).$

\begin{prop}\label{formenormale}Let $\mu$ be any non zero number. Every 3-web $(x,y,f),$ with non-zero curvature at the origin, is locally isomorphic to an unique 3-web $(x,y,F_{\mu})$ such that 
$$F_\mu (x,y)=x+y+xy(x-y)(\mu+g(x,y))$$
 where $g$ is a function vanishing at the origin.
\end{prop}

This proposition is a particular case of the existence and unicity of {\sl normal form} for 3-webs which appears for the first time in
\cite{DJ}: for any 3-web $V$ near any point $m,$ there are local coordinates $(x,y)$ vanishing at $m$ such that $V$ becomes $(x,y,h)$  with
$$h(x,y)=x+y+xy(x-y)k(x,y)$$
where $k$  may be any function. Moreover $k$ is unique up to a homothety $(x,y)\mapsto (\lambda x,\lambda y).$ For the moment we impose no assumption concerning curvature. If the curvature doesn't vanish at $m$ then $k(0,0)$ is different from zero. Then, up to a homothety, we can assume $k(0,0)=\mu$ and we obtain the above proposition. 

The fact that $W$ and $\bar W$ are non flat implies that $(x,y,f_W)$ and $(x,y,f_{\bar W})$ have non zero curvature at the origin. The above proposition says that there are two  functions $ g_{\mu}$ and $\bar g_{\mu}$, vanishing at the origin, such that $W$ and $\bar W$ are respectively isomorphics to $(x,y,F_{\mu})$ and $(x,y,\bar F_{\mu}),$ with
$$F_{\mu}(x,y)=x+y+xy(x-y)(\mu+g_{\mu}(x,y)),$$
$$\bar F_{\mu}(x,y)=x+y+xy(x-y)(\mu+\bar g_{\mu}(x,y)).$$

The unicity part of the above proposition implies the following lemma.
\begin{lemma} Fix the non zero number $\mu .$ Then $W$ is isomophic to $\bar W,$ by an isomorphism which maps $M$ to $\bar M,$ if and only if 
$$ g_{\mu} =\bar g_{\mu}.$$ \end{lemma}

To be able to express the isomorphy of $W$ and $\bar W$  with this lemma, we need to have a practical method  to compute  $g_{\mu}$ and $\bar g_{\mu}.$

We don't have such a method but, at least, we will give in the next section an algorithm which  computes the $k$-jet of $g_{\mu}$ (resp. $\bar g_{\mu}$), starting with $(k+2)$-jets of $a$, $b$ and $c$ (resp. $\bar a$, $\bar b$ and $\bar c$).   If $W$ and $\bar W$ are isomorphic and for every $k$, this will gives many precise polynomial equations between  $(a_2,b_2,c_2,\cdots ,a_{k+2}, b_{k+2}, c_{k+2})$ and
 $({\bar a}_2,{\bar b}_2,{\bar c}_2,\cdots,{\bar a}_{k+2},{\bar b}_{k+2},{\bar c}_{k+2}).$

\section{Algorithm to compute $k$-jets of $g_{\mu}.$}

In this section we mimic a classical proof of the existence of {\sl normal form} for any 3-web to get an algorithm which works for jets.

We choose 
$$\mu=a_2+b_2+c_2$$
i.e. $\mu$  is the characteristic of $W$ at the origin up to the factor 4.

The {\bf input} is the $(k+2)$-jet of $a$, $b$ and $c$, i.e.  the numbers $a_2,\cdots ,a_{(k+2)},$  $b_2,\cdots ,b_{(k+2)}$ and  $c_2,\cdots ,c_{(k+2)}.$ 

The {\bf first procedure} gives the $(k+3)$-jet, $j^{(k+3)}f_W$ of $f_W$ at the origin by computing the $(k+3)$-jet of the solution of the implicit relation 
(\ref{pointcommun}) when we replace respectively $a,$ $b$ and $c$ by their $(k+2)$-jets.

The {\bf second procedure} is the {\sl normalisation} of $j^{(k+3)}f_W$ to obtain the $k$-jet of $g_\mu.$ We do that in six  steps.

{\sl First step}. We use the simplifying notation $F=j^{(k+3)}f_W$ and consider it as a polynomial function with two variables.  Compute $X,$ the $(k+3)$-jet of the inverse function of $t\mapsto F(t,0).$

{\sl Second step}.  Compute $Y,$ the $(k+3)$-jet of the inverse function of $t\mapsto F(0,t).$

{\sl Third step}. Compute the $(k+3)$-jet of $F(X(x),Y(y)).$ We denote it by $G$ (remark that $G$ has the shape $x+y+xy\Theta (x,y)$)

{\sl Fourth step}. Consider $K=G(t,t).$ Find the 1-variable polynomial $U=t+u_2t^2+\cdots  u_{(k+3^)}t^{(k+3)}$ such that
the $(k+3)$-jet of $G(U(t),U(t))$ is egal to $U(2t)$ (it exists by the classical Sternberg's theorem \cite{SS} which says that every map $t\mapsto 2t+d_2t^2+\cdots$ is conjugated to its linear part ; it is also unique because $j^1U=t$).

{\sl Fifth step}. Compute $V,$ the $(k+3)$-jet of the inverse of $U,$ and $H,$  the $(k+3)$-jet of $V(G(U(x),U(y)))$ (remark that $H$ has the shape $x+y+xy(x-y)\Psi (x,y)$).

{\sl Sixth step}. Compute $ L,$ the $k$-jet of $(H-x-y)/(xy(x-y))$.

The {\bf output} is $E=L-\mu$ which is the k-jet of $g_{\mu}.$

We have implemented this algorithm on Maple. It works very rapidly if $k$ is less or equal to seven.

\section{Proof of Theorem \ref{resultatA}.}
We keep notation $\mu:=a_2+b_2+c_2.$
The assumption of our theorem writes as 
$$\bar a_2+\bar b_2+\bar c_2=\mu.$$

Using the algorithm of the preceeding section we compute respectively the 5-jets $E$  of  $g_{\mu}$ and   $\bar E$ of $\bar
g_{\mu}.$

We write the Taylor expansion of $E$  as $E_{10} x+ E_{01}y+\cdots+E_{ij}x^iy^j.$ We have
\begin{equation}\label{rel1} E_{10}=(-2a_2-2b_2+c_2+20a_3+8b_3+14c_3)/7.\end{equation}

\begin{equation}\label{rel2}E_{01}=(2a_2+2b_2-c_2+20b_3+8a_3+14c_3)/7.\end{equation}

And also

$$ E_{20}=a_2/12+b_2/12-c_2/6+10a_2^2/3+c_2^2+2a_2b_2/3-2b_2c_2/3+2a_2c_2/3$$ $$-2a_3+c_3/3+20a_4/3+4b_4/3+3c_4.$$

$$E_{11}=-a_2/12-b_2/12-c_2/3+a_2^2/3-b_2^2/3+4c_2(a_2-b_2)/3+4a_4+4b_4+5c_4.$$

$$ E_{02}=a_2/12+b_2/12-c_2/6-10b_2^2/3-c_2^2-2a_2b_2/3-2b_2c_2/3+2a_2c_2/3$$ $$+2b_3-c_3/3+20b_4/3+4a_4/3+3c_4.$$

The following $ E_{ij}$ may  have very long expressions. For example if $i+j =5$ they contain nearly hundred terms. We only retain that they are polynomial expressions, with rational coefficients, in some of the $a_r,$ $b_r$ and $c_r$ variables.

In the case of $\bar E$, we obtain the same expressions for its coefficients ${\bar E}_{ij}$ except we have to change respectivily $a_r,$ $b_r$ and $c_r$  by  $\bar a_r,$ $\bar b_r$ and $\bar c_r.$ 

To simplify  we use the following notations:
$$\bar a_r=a_r+A_r,\ \ \bar b_r=b_r+B_r,\ \  \bar c_r=c_r+C_r,$$ 
for every $r.$
The hypothesis of our theorem writes as $C_2=-A_2-B_2.$ To prove it we have to prove  $A_r=B_r=C_r=0$ for every $r.$

We adopt notations
$$T_{ij}=E_{ij}-{\bar E}_{ij},$$
for every $i$ and $j.$ The existence of an isomorphism between $W$ and $\bar W$ implies the set of equations $T_{ij}=0.$ They are polynomial equations with unknown  $A_r,$ $B_r$ and $C_r$ and coefficients rational in some of the $a_r,$ $b_r,$ $c_r.$

For example the relations (\ref{rel1}) and  (\ref{rel2}) give $T_{10}$ and $T_{01},$ i.e. the order one equations. Equations 
$T_{10}=0$ and $T_{01}=0,$ give relations
$$A_3=A_2/4+B_2/4-C_3/2,$$
$$B_3=-A_2/4-B_2/4-C_3/2.$$

At order 2 the equations $T_{20}=0,$ $T_{11}=0$ and $T_{02}=0$ give
$$A_4 = A_2/8+c_2B_2/3+b_2B_2/3+B_2/8-b_2A_2/12-B_2A_2/2-B_2a_2/6$$$$-19a_2A_2/12-A_2^2+5c_2A_2/12-C_3/4,$$
$$B_4 = A_2/8-5c_2B_2/12+19b_2B_2/12+B_2/8-a_2A_2/3+b_2A_2/6$$$$+B_2A_2/2+B_2a_2/12+B_2^2-c_2A2/3+C_3/4,$$
$$C_4 = -A_2/4+c_2B_2/3-5b_2B_2/3-B_2/4+5a_2A_2/3-b_2A_2/3$$$$+B_2a_2/3-B_2^2-c_2A_2/3+A_2^2.$$
Note that the second members of these relations are polynomial with variables $A_2,$ $B_2,$ $C_3,$ $a_2,$ $b_2$ and $c_2.$

At order 3 and using Maple, we obtain four equations $T_{ij}=0$ and they allow to express $A_5,$ $B_5,$ $C_5$ and $C_3$ as rational fonctions of the variables $A_2,$ $B_2,$ $a_2,$ $b_2,$ $c_2,$ $a_3,$ $b_3$ and $c_3,$ with the denominator $a_2+b_2+c_2.$ Formulas are two long to be reproduced here.

For the moment we skip order 4 and consider the six equations $T_{ij}=0$ with $i+j=5.$  Maple proves that they allow to express $A_6,$ $B_6,$ $C_6,$  $A_7,$ $B_7$ and $C_7$ rationally in function of $A_2,$ $B_2,$ $a_2,$ $b_2,$ $c_2,$ $\cdots ,$ $a_5,$ $b_5$ and $c_5$ with the denominator $(a_2+b_2+c_2)^2.$

We remark that $A_2=B_2=0$ implies that $A_i,$ $B_i$ and $C_i$ vanish for $i=3,\cdots ,7.$

Now we compute the five equations $T_{ij}=0$ with $i+j=4.$ They have the shape
$$\alpha_j A_2^2+\beta_j B_2^2+\gamma_j A_2B_2+\mu_j A_2+\nu_j B_2 =0$$
where $\alpha_j,\beta_j ,\gamma_j ,\mu_j ,\nu_j $ are rationally functions of $b_2,$ $c_2,$ $\cdots ,$ $a_5,$ $b_5$ and $c_5$ with the denominator $a_2+b_2+c_2.$

We see also that the $(3\times 5)$-matrix with lines $(\alpha_j,\beta_j ,\gamma_j )$ has rank 3 for any value of $b_2,$ $c_2,$ $\cdots ,$ $a_5,$ $b_5$ and $c_5.$  This is a consequence of the fact that the coefficients of this matrix depends only on the three numbers $a_2+b_2+c_2,$ $a_3+b_3+c_3$ and $a_2+b_2+2b_3-2a_3.$
So our system of equations can be rewriten as
\begin{equation}\label{ordre4}\matrix{A_2^2=\psi_1A_2+\phi_1B_2,\cr
B_2^2=\psi_2A_2+\phi_2B_2,\cr
A_2B_2=\psi_3A_2+\phi_3B_2,\cr
0=\psi_4A_2+\phi_4B_2,\cr
0=\psi_5A_2+\phi_5B_2. }\end{equation}

{\bf To obtain the Maple worksheet which gives these results, contact the author at dufourh@netcourrier.com}

\begin{lemma} For $i+j>2$ we have relations 
$$E_ {ij}=\theta_{ij}a_{i+j+2}+\phi_{ij}b_{i+j+2}+\psi_{ij}c_{i+j+2}+S_{ij}$$ where $\theta_{ij},$ $\phi_{ij}$ and $\psi_{ij}$ are some constants and $S_{ij}$ is a polynomial with variables $a_2,$  $b_2,$ $c_2,$ $\cdots ,$ $a_{i+j+1},$ $b_{i+j+1}$ and $c_{i+j+1}.$ Moreover any $(3\times 3)$-submatrix of the matrix whose lines are $(\theta_{ij},\phi_{ij},\psi_{ij})$ is of rank 3.

\end{lemma}

This can be proven for any $(i,j)$ (whithout Maple !) as follows. We use notation $n=i+j+2.$  We first see that $F,$ the $(n+1)$-jet of $f_W,$ has the shape
$$(x+y)/2+ \Theta +a_nK+b_nL+c_nM,$$
where $K,$ $L$ and $M$ are homogeneous polynomials of degree $(n+1),$ with variable $(x,y)$ and constant coefficients;
$\Theta$ is a polynomial expression with variables $x,$ $y,$  $a_2,$  $b_2,$ $c_2,$ $\cdots ,$ $a_{n-1},$ $b_{n-1}$  and $c_{n-1}.$ Now we apply  the normalising procedure, described in the previous section, to $F.$ We only have to follow what happens to the terms containing $a_n,$ $b_n$ and $c_n.$ It is a little long but elementary.

 Note that, with this lemma, we recover in part  above  Maple results: we recover that  the equations $T_{ij}=0$ for $i+j=n-2$ 
give expressions of $A_n,$ $B_n$ and $C_n$ in terms of the previous $A_p,$ $B_p$ and $C_p,$ and some of the $a_r$, $b_r$ and $c_r,$ for $n=4,5,6,7.$

This lemma proves also that, if $A_2,$ $B_2,$ $A_3,$ $B_3$ and $C_3$ vanish then all the $A_r,$ $B_r$ and $C_r$ vanish also. Using  Maple calculations above we see that the relation $A_2=B_2=0$ implies that all the $A_r,$ $B_r$ and $C_r$ vanish but also all the $T_{ij}.$

Using the three first equations (5) we can replace any monomial in $A_2$ and $B_2$ by a linear expression $\rho A_2+\tau B_2.$ So our system of equations $T_{ij}$ can rewriten as a set of equations

$$A_r=u_r^1A_2+u_r^2B_2,\ B_r=v_r^1A_2+v_r^2B_2,\ C_r=w_r^1A_2+w_r^2B_2,\ t^1_kA_2+t^2_kB_2=0;$$
for $ r>2$ and an infinity of $k;$ the coefficients $u^s_r,$ $v^s_r,$ $w^s_r$ and $t^s_r$ depending only on the $a_n,$ $b_n$ and $c_n.$

If this system has  a non zero solution $(A_2,B_2,A_3,B_3,C_3,\dots)$ then it has an infinity of solutions: $t(A_2,B_2,A_3,B_3,C_3,\dots)$ for any number $t.$

This means that $W$ would be isomorphic to any linear 3-web $\bar W_t$ which is described by the three one variable functions
$$\bar a_t=a+tA,\ \bar b_t=b+tB,\ \bar c_t=c+tA,$$
where $A,$ $B$ and $C$ are the functions which have respectively Taylor expansions $A_2x^2+\cdots +A_nx^n+\cdots,$
 $B_2y^2+\cdots +b_ny^n+\cdots$ and  $C_2z^2+\cdots +C_nz^n+\cdots .$

This contredicts the known fact that any non flat linear 3-web can only be isomorphic to a finite number of  homographically different linear 3-webs.

So we have proven that the only possiblity is $A_r=B_r=C_r=0$ for any $r$ and  our Theorem.

\section{Remarks.}

Let $\phi$ be an isomorphism between the two linear non flat webs $W$ and $\bar W.$ We suppose also that $\phi$ maps a point $M$ of the domain of $W$ to a point $\bar M$ on the domain of $\bar W.$

We adopt the description of $W$ (resp. $\bar W$) near $M$ (resp. $\bar M$) of the section \ref{descro}, i.e. with three 1-variable functions $a,$ $b$ and $c$ (resp.  $\bar a,$ $\bar b$ and $\bar c$). Then $\phi$ becomes a local diffeomorphism which preserves the origin. As it preserves the $u$-axis and the two bissectrices $u=v$ , $u=-v,$ its 1-jet at the origin is a homothety $kI$. Then we have
$$car_W(M)=k^2car_{\bar W}(\bar M).$$

For any linear 3-web we find in \cite{SA} the construction of 1-forms $U_1,$ $U_2$ and $U_3 $ which are invariant up to homographies, such that the three foliations  are given by the kernels of these forms and 
$$U_1+U_2+U_3=0.$$ We denote  $U_1,$ $U_2$ and $U_3 $ these forms for $W$ and  $\bar U_1,$ $\bar U_2$ and $\bar U_3 $
for $\bar W.$
  
Classically, there is a function $f$ such that
$$\phi^{*}\bar U_i=f.U_i$$
for every $i=1,2,3.$ Using the  description of section \ref{descro}, we can show 
$$car_W(M)=f^2car_{\bar W}(\bar M).$$

So our result can be rephrased as: if $f^2$ is equal to 1 at some point then, near this point, $\phi$ is a homography.

\end{document}